\documentclass[dvipsnames,10pt]{article}
\usepackage{graphicx,float,soul}
\usepackage{amsmath, amsfonts, amsthm, amssymb}
\usepackage{braket, nicefrac}
\usepackage{hyperref,mathtools}
\usepackage{tikz}
\usepackage{tikz-cd}
\usepackage{theoremref, verbatim}
\usepackage[a4paper, total={5.5in, 9.1in}]{geometry}
\usepackage{setspace}
\usepackage{xcolor}
\usepackage{tgpagella}
\usepackage{needspace}

\newcommand{\lcm}{\mathrm{lcm}}

\newcommand{\reg}{\mathrm{reg}}
\newtheorem{theorem}{Theorem}[section]
\newtheorem{lemma}[theorem]{Lemma}
\newtheorem{proposition}[theorem]{Proposition}
\newtheorem{cor}[theorem]{Corollary}
\newtheorem{defi}[theorem]{Definition}
\newtheorem{construction}[theorem]{Construction}
\newtheorem{example}[theorem]{Example}
\newtheorem{remark}[theorem]{Remark}
\title{\Large \textbf{Generalized Edge Ideals and Regularity Offset of Graphs}}

\author{
    Grigoriy Blekherman \thanks{Department of Mathematics, Georgia Institute of Technology, greg@math.gatech.edu} 
    \and 
    Jiaheng Li\thanks{Department of Mathematics, Texas A\&M University, jli3440@tamu.edu} 
    \and 
    Jianuo Zhou\thanks{Department of Mathematics, Georgia Institute of Technology, jianuozhou2@outlook.com} 
}

\newcommand{\mycomment}[1]{}

\date{}

\begin{document}

\setstretch{1.2}
\maketitle

\begin{abstract}
We introduce generalized edge ideals by assigning pairwise coprime monomials, instead of variables, to the vertices of a graph. We use the lcm-lattice to show that the multigraded Betti numbers and Castelnuovo-Mumford regularity of the quotient by a generalized edge ideal agree with the Betti numbers and resolution regularity of the quotient by an ordinary edge ideal in a non-standard $\mathbb{Z}$-grading. We define the regularity offset as a measurement of the difference between the regularities of the quotients by a generalized edge ideal and its ordinary counterpart, and prove bounds for the regularity offset. We compute the regularity offset for co-chordal graphs and paths, and derive a formula for the regularity offset of forests. We also derive a graph-theoretic criterion for when a forest's regularity offset is always zero, or equivalently, its top-level Betti number is nonzero.
\end{abstract}



\section{Introduction}

Quadratic square-free monomial ideals correspond bijectively to simple graphs. The vertices are labeled by variables, and the ideal is generated by monomials corresponding to the graph's edges (or, under an alternative convention, its non-edges). Algebraic properties of edge ideals, such as the minimal free resolution, Betti numbers, and Castelnuovo-Mumford regularity, have been studied extensively in \cite{blekherman2020bounds,bouchat,dao2013bounds,fakhari2017regularity,ha2014regularity,ha2008monomial,kozlov1999complexes,morales2012regularity,woodroofe2014matchings}, to name just a few of the papers. 

We introduce generalized edge ideals, where the vertices of a graph $G$ are labeled by pairwise coprime monomials instead of single variables, and the ideal is generated by products of vertex monomials corresponding to the edges of $G$. Our main interest is to see how various algebraic parameters of the quotients by generalized and ordinary edge ideals are related. We use the regularity conventions in Section \ref{notation}. Our first crucial observation, using the lcm-lattice introduced in \cite{welker1999lcm}, is that the minimal free resolutions and graded Betti numbers of the quotients by generalized edge ideals correspond to those of ordinary edge ideal quotients where the base polynomial ring is given a positive $\mathbb{Z}$-grading. Consequently, their regularities agree under the Betti-degree definition. We illustrate this with an example:

\Needspace{10\baselineskip}
\begin{example}\thlabel{firstexample}
Consider the following graph:
\begin{center}

\tikzset{every picture/.style={line width=0.75pt}} 

\begin{tikzpicture}[x=0.75pt,y=0.75pt,yscale=-1,xscale=1]

\draw    (191,70.33) -- (230.1,91.33) ;
\draw  [fill={rgb, 255:red, 0; green, 0; blue, 0 }  ,fill opacity=1 ] (188.45,70.33) .. controls (188.45,68.93) and (189.59,67.79) .. (191,67.79) .. controls (192.41,67.79) and (193.55,68.93) .. (193.55,70.33) .. controls (193.55,71.74) and (192.41,72.88) .. (191,72.88) .. controls (189.59,72.88) and (188.45,71.74) .. (188.45,70.33) -- cycle ;
\draw  [fill={rgb, 255:red, 0; green, 0; blue, 0 }  ,fill opacity=1 ] (227.55,91.33) .. controls (227.55,89.93) and (228.69,88.79) .. (230.1,88.79) .. controls (231.5,88.79) and (232.64,89.93) .. (232.64,91.33) .. controls (232.64,92.74) and (231.5,93.88) .. (230.1,93.88) .. controls (228.69,93.88) and (227.55,92.74) .. (227.55,91.33) -- cycle ;
\draw  [fill={rgb, 255:red, 0; green, 0; blue, 0 }  ,fill opacity=1 ] (227.55,50) .. controls (227.55,48.59) and (228.69,47.45) .. (230.1,47.45) .. controls (231.5,47.45) and (232.64,48.59) .. (232.64,50) .. controls (232.64,51.41) and (231.5,52.55) .. (230.1,52.55) .. controls (228.69,52.55) and (227.55,51.41) .. (227.55,50) -- cycle ;
\draw    (230.1,50) -- (191,70.33) ;
\draw    (269.19,71) -- (230.1,91.33) ;
\draw    (230.1,50) -- (269.19,71) ;
\draw  [fill={rgb, 255:red, 0; green, 0; blue, 0 }  ,fill opacity=1 ] (266.64,71) .. controls (266.64,69.59) and (267.78,68.45) .. (269.19,68.45) .. controls (270.6,68.45) and (271.74,69.59) .. (271.74,71) .. controls (271.74,72.41) and (270.6,73.55) .. (269.19,73.55) .. controls (267.78,73.55) and (266.64,72.41) .. (266.64,71) -- cycle ;
\draw    (230.1,50) -- (230.1,91.33) ;
\draw    (345.67,70.33) -- (384.76,91.33) ;
\draw  [fill={rgb, 255:red, 0; green, 0; blue, 0 }  ,fill opacity=1 ] (343.12,70.33) .. controls (343.12,68.93) and (344.26,67.79) .. (345.67,67.79) .. controls (347.07,67.79) and (348.21,68.93) .. (348.21,70.33) .. controls (348.21,71.74) and (347.07,72.88) .. (345.67,72.88) .. controls (344.26,72.88) and (343.12,71.74) .. (343.12,70.33) -- cycle ;
\draw  [fill={rgb, 255:red, 0; green, 0; blue, 0 }  ,fill opacity=1 ] (382.21,91.33) .. controls (382.21,89.93) and (383.35,88.79) .. (384.76,88.79) .. controls (386.17,88.79) and (387.31,89.93) .. (387.31,91.33) .. controls (387.31,92.74) and (386.17,93.88) .. (384.76,93.88) .. controls (383.35,93.88) and (382.21,92.74) .. (382.21,91.33) -- cycle ;
\draw  [fill={rgb, 255:red, 0; green, 0; blue, 0 }  ,fill opacity=1 ] (382.21,50) .. controls (382.21,48.59) and (383.35,47.45) .. (384.76,47.45) .. controls (386.17,47.45) and (387.31,48.59) .. (387.31,50) .. controls (387.31,51.41) and (386.17,52.55) .. (384.76,52.55) .. controls (383.35,52.55) and (382.21,51.41) .. (382.21,50) -- cycle ;
\draw    (384.76,50) -- (345.67,70.33) ;
\draw    (423.86,71) -- (384.76,91.33) ;
\draw    (384.76,50) -- (423.86,71) ;
\draw  [fill={rgb, 255:red, 0; green, 0; blue, 0 }  ,fill opacity=1 ] (421.31,71) .. controls (421.31,69.59) and (422.45,68.45) .. (423.86,68.45) .. controls (425.26,68.45) and (426.4,69.59) .. (426.4,71) .. controls (426.4,72.41) and (425.26,73.55) .. (423.86,73.55) .. controls (422.45,73.55) and (421.31,72.41) .. (421.31,71) -- cycle ;
\draw    (384.76,50) -- (384.76,91.33) ;

\draw (167.81,66.23) node [anchor=north west][inner sep=0.75pt]    {$x_{1}$};
\draw (222.48,101.07) node [anchor=north west][inner sep=0.75pt]    {$x_{2}$};
\draw (221.81,32.57) node [anchor=north west][inner sep=0.75pt]    {$x_{3}$};
\draw (276.48,66.9) node [anchor=north west][inner sep=0.75pt]    {$x_{4}$};
\draw (300.48,79.23) node [anchor=north west][inner sep=0.75pt]    {$m_1=z_{1} z_{2}$};
\draw (377.48,31.9) node [anchor=north west][inner sep=0.75pt]    {$m_3=z_{3}$};
\draw (371.14,98.4) node [anchor=north west][inner sep=0.75pt]    {$m_2=z_{4} z_{5}^{2}$};
\draw (418.81,74.23) node [anchor=north west][inner sep=0.75pt]    {$m_4=z_{6}^{3}$};
\end{tikzpicture}

\end{center}
\begin{equation*}
    I=( x_{1} x_{2} ,x_{1} x_{3} ,x_{2} x_{3} ,x_{2} x_{4} ,x_{3} x_{4})
\end{equation*}
\begin{equation*}
    J=\left( z_{1} z_{2} z_{3} ,z_{1} z_{2} z_{4} z_{5}^{2} ,z_{3} z_{4} z_{5}^{2} ,z_{3} z_{6}^{3} ,z_{4} z_{5}^{2} z_{6}^{3}\right), \,\,\,\,  {\deg z_{i} =1}
\end{equation*}
\begin{equation*}
    J'=( y_{1} y_{2} ,y_{1} y_{3} ,y_{2} y_{3} ,y_{2} y_{4} ,y_{3} y_{4}), \,\,\, \deg y_{1} =2,\deg y_{2} =3,\deg y_{3} =1,\deg y_{4} =3
\end{equation*}

Labeling vertices with distinct variables gives the ordinary edge ideal $I$. Labeling the vertices by pairwise coprime monomials $m_i$ defines a generalized edge ideal $J$ in a larger polynomial ring.

Construct an ideal $J'$ that lies in the ring $k[y_1,y_2,y_3,y_4]$, which is equipped with the grading $\deg y_i=\deg m_i$. The new ideal $J'$ is an ordinary edge ideal in a ring with non-standard $\mathbb{Z}$-grading. The quotients $k[y_1,y_2,y_3,y_4]/J'$ and $k[z_1,\cdots,z_6]/J$ have corresponding minimal free resolutions, the same graded Betti numbers, and the same regularity in the sense of Section \ref{notation}.

\end{example}

We then concentrate on the relationship between the regularities of the quotients by ordinary and generalized edge ideals. One simple but important insight is an inequality which bounds their difference.

\begin{theorem}[\thref{plusn}]\thlabel{thm:offset}
    Let $G$ be a graph with $V(G)=\{x_1,\cdots,x_n\}$. Let $I_G$ be the ordinary edge ideal of a graph $G$ and $J_G$ be its edge ideal under the grading $\mathbf{g}=(k_1,\cdots,k_n)$. Let
    \begin{equation*}
        N=\sum_{i=1}^n(k_i-1)
    \end{equation*}
    be the total increment of the degrees of all variables. Then,
    \begin{equation*}
        \mathrm{reg}(G) \leqslant \mathrm{reg}(G,\mathbf{g}) \leqslant \mathrm{reg}(G)+N.
    \end{equation*}
\end{theorem}

In many cases the upper bound in the above theorem is tight, i.e. $\mathrm{reg}(G,\mathbf{g})$ is determined by $\mathrm{reg}(G)$ and the total added degree in the new grading. To study this phenomenon we introduce the \emph{regularity offset} of a graph under a grading:


\begin{defi}[Regularity Offset, \thref{ro}]
    Continue with the notation of the previous theorem. We define the regularity offset of $G$ under grading $\mathbf{g}$ as
    \begin{equation*}
        \mathrm{ro}(G,\mathbf{g})=\mathrm{reg}(G)+N-\mathrm{reg}(G,\mathbf{g})
    \end{equation*}
    We say that $G$ does not have regularity offset if $\mathrm{ro}(G,\mathbf{g})= 0$ for all gradings $\mathbf{g}\in (\mathbb{N}^+)^n$. In other words, $\mathrm{reg}(G,\mathbf{g})$ always attains maximum growth.
\end{defi}

\begin{example}
Let $G$ be the graph in \thref{firstexample}, let $\mathbf{g}=(2,3,1,3)$, and write $R=k[z_1,\cdots,z_6]$ with its standard grading. Then $\mathrm{reg}(G)=1$ and $\mathrm{reg}(R/J)=\mathrm{reg}(G,\mathbf{g})=6$. Therefore, the regularity offset is $1+5-6=0$.
\end{example}

\begin{remark}
We can define the operation of substituting a variable with a product of disjoint variables for any monomial ideal. Understanding what this does to regularity reduces to examining positive $\mathbb{Z}$-gradings in exactly the same way, and a more general version of \thref{thm:offset} also holds with the same proof. If $x_1^{c_1}\cdots x_n^{c_n}$ is the lcm of all generators, the bound would be
\begin{equation*}
    \mathrm{reg}_{\mathrm{g}}(S/I) \leqslant \mathrm{reg}(S/I) + \sum_{i=1}^n c_i(k_i-1).
\end{equation*}
Therefore, the notion of regularity offset also makes sense for arbitrary monomial ideals. In this paper we only analyze edge ideals, as this allows us to take advantage of graph structure.
\end{remark}

As indicated in the definition of regularity offset, some graphs will have maximal growth of regularity \emph{regardless of the grading}. These are precisely the graphs with no regularity offset. The key to understanding regularity offset is the following Lemma which shows that to have no regularity offset is equivalent to maximal growth of regularity for a \emph{single grading} where each degree is increased by one.

\begin{lemma}[\thref{rogt0}]
    A graph $G$ with $n$ vertices does not have regularity offset if and only if
    \begin{equation*}
        \mathrm{ro}(G,\mathbf{t})=0,
    \end{equation*}
     where $\mathbf{t}\in\mathbb{Z}^n$, $\mathbf{t}=(2,\cdots,2)$.
\end{lemma}
Together with Hochster's formula, this allows us to give a topological characterization of graphs with no regularity offset in \thref{wholegraph}. We also show in \thref{prlower} that there exist graphs and gradings for which the \emph{lower bound} of \thref{plusn} is attained, for any number of added degrees.

In Section \ref{rossg} we use the tools above to calculate regularity offset of co-chordal graphs, paths and cycles.
A co-chordal graph is a graph whose complement is a chordal graph. 

\begin{theorem}[\thref{ccd}]
    Let $G$ be co-chordal, let $\mathbf{g}$ be a positive grading, and put $t_v=\deg_{\mathbf{g}}v-1$. Let $\mathcal{C}_{\min}(G)$ consist of the inclusion-minimal sets $C\subseteq V(G)$ for which $G^c[V(G)\setminus C]$ has at least two connected components. Then
    \[
        \mathrm{ro}(G,\mathbf{g})=
        \begin{cases}
            \displaystyle\min_{C\in\mathcal{C}_{\min}(G)}\sum_{v\in C}t_v,
                & E(G)\neq\varnothing,\\[4pt]
            \displaystyle\sum_{v\in V(G)}t_v,
                & E(G)=\varnothing.
        \end{cases}
    \]
    If $G^c$ has at least two connected components, the empty set is the unique member of $\mathcal{C}_{\min}(G)$, so the offset is zero.
\end{theorem}

Using results of \cite{kozlov1999complexes}, we determine the behavior of regularity offset for paths and cycles.

\begin{proposition}[\thref{cycle}, \thref{path}] \thlabel{prop17}
    The cycle of length $n$ does not have regularity offset. The path with $n$ vertices has regularity offset if and only if $n\equiv 1$ $(\mathrm{mod}$ $3)$.
\end{proposition}
For paths with non-zero regularity offset we provide more details about the offset of individual gradings.
\begin{theorem}[\thref{path2}] \thlabel{thm18}
    Let $n=3k+1$ with $k\geq0$, and label the vertices of $P_n$ in order as $x_1,\ldots,x_n$. For a positive grading $\mathbf{g}$, put $t_i=\deg_{\mathbf{g}}x_i-1$. Then $\mathrm{ro}(P_n,\mathbf{g})=\min_{0\leq s\leq k}t_{3s+1}$.
\end{theorem}

In the final section we move on to the regularity offset of forests. We first obtain a graph-theoretic criterion for when a forest does not have regularity offset.

\begin{proposition}[\thref{deletion2}, \thref{deletion3}]
    Define a deletion operation as: take a leaf with a neighbor, and remove all vertices in the neighborhood of its neighbor, plus the neighbor itself.

    Let $T$ be a forest. It has no regularity offset if and only if it can be reduced to an empty tree after a finite number of deletion operations. We also call these forests ``full forests''.

    When $T$ is full, there exists a unique $t$ such that $\beta_{t,\mathbf{1}}(I_T)\neq 0$. Denote this number by $i_T$. Then $i_T=|V(T)|-\reg(T)$.
\end{proposition}

\begin{remark}
    When applied to paths, this criterion recovers the path statement of \thref{prop17}.
\end{remark}

This leads us to the final result, on the regularity of forests.

\begin{theorem}[\thref{final}]
    Let $T$ be a forest on $\{x_1,\cdots,x_n\}$ and $Z$ the set of full induced subforests of $T$. For a grading $\mathbf{g}$, let $t_i=\deg x_i-1$ be the increment of the degree of each $x_i$. Then we have
    \begin{equation*}
        \mathrm{ro}(T,\mathbf{g}) = \mathrm{reg}(T) + \min_{S\in Z\;\text{maximal}}\left( \sum_{i:\,x_i\notin V(S)} t_i - \mathrm{reg}(S) \right).
    \end{equation*}
\end{theorem}

\begin{remark}
    When applied to paths, this formula recovers \thref{thm18}.
\end{remark}

\section{Preliminaries}

\subsection{Notation}\label{notation}

Throughout the paper, $k$ denotes a field of characteristic $0$. We use both $\mathbb{Z}$-grading and $\mathbb{Z}^n$-grading (= multigrading) for the base rings of edge ideals. For the ring $S=k[x_1,\cdots,x_n]$, we specify gradings in the form of $\mathbf{g}=(k_1,\cdots,k_n)\in (\mathbb{Z}^+)^n$. When using $\mathbb{Z}$-grading, the grading $\mathbf{g}$ specifies that the degree of $x_i$ is $k_i$. When using $\mathbb{Z}^n$-grading, $\mathbf{g}$ specifies that the degree of $x_i$ is $k_i\mathbf{e}_i$, where $\mathbf{e}_i$ is the $i$-th unit vector. We will always say $\deg x_i=k_i$, even when using the $\mathbb{Z}^n$-grading, as this does not lead to any ambiguity.

We use $G$ to denote a graph, usually with $n$ vertices $x_1,\cdots,x_n$. We also use $T$ for forests and $P_n$ for the path with $n$ vertices. We use $I_G$ to denote the edge ideal of $G$, whose ambient ring is $k[G]:= k[V(G)]$. With the standard grading on $k[G]$, we write $\mathrm{reg}(G)=\mathrm{reg}(k[G]/I_G)$ for the Castelnuovo-Mumford regularity of the quotient ring. Thus $\mathrm{reg}(G)=0$ when $G$ has no edges, including when $G$ is empty.

We denote the edge ideal of $G$ under a positive grading $\mathbf{g}$ by $J_G$. For the quotient ring equipped with this grading, we set
\[
    \mathrm{reg}(G,\mathbf{g})=\mathrm{reg}(k[G]/J_G)
    :=\max\{j-i\mid\beta_{i,j}(k[G]/J_G)\neq0\}.
\]
This is resolution regularity. In the standard grading it agrees with Castelnuovo-Mumford regularity, so $\mathrm{reg}(G,\mathbf{1})=\mathrm{reg}(G)$. For other positive gradings, this convention need not agree with the definition using local cohomology.

Our Betti-number notation $\beta_{i,\mathbf{a}}(I)$ and $\beta_{i,j}(I)$ refers to the quotient $S/I$, as in the resolution in \thref{mfr}. Every regularity value in this paper likewise refers to an explicitly specified quotient ring or to the graph notation just defined.

\subsection{The minimal free resolution and Betti numbers}

\begin{defi}\thlabel{mfr}
    The \textbf{minimal free resolution} of a homogeneous ideal $I\subset S=k[x_1,\cdots,x_n]$ is a finite exact sequence of graded maps
    \begin{equation*}
        \begin{tikzcd}
            0 \arrow[r, ""] & \mathbf{F}_m \arrow[r, "{d_m}"] & \mathbf{F}_{m-1} \arrow[r, ""] & \cdots \arrow[r, ""] & \mathbf{F}_1 \arrow[r, "{d_1}"] & (\mathbf{F}_0 = S) \arrow[r, "{d_0}", two heads] & S/I \arrow[r, ""] & 0
        \end{tikzcd}
    \end{equation*}
    where each $\mathbf{F}_i$ is a free $S$-module, and $d_i(\mathbf{F}_i)\subset (x_1,\cdots,x_n)\mathbf{F}_{i-1}$ for all $1\leqslant i\leqslant m$. Equivalently, let $M(d_i)$ be the matrix representation of $d_i$ (so $M(d_i)$ is a matrix of elements in $S$), then every entry in every $M(d_i)$ must lie in the ideal $(x_1,\cdots,x_n)$.
\end{defi}

By Hilbert's syzygy theorem \cite{hilbertsyzygy}, minimal free resolutions of ideals in polynomial rings always exist, and are unique up to isomorphism. When $S$ is graded, the maps $d_i$ automatically become graded maps. These degrees of the generators of each module $\mathbf{F}_i$ are recorded by the Betti numbers.

\begin{defi}
    Let $S=k[x_1,\cdots,x_n]$ with grading $\mathbf{g}=(k_1,\cdots,k_n)$. We use the notation of \thref{mfr}. Let $\mathbf{a}\in\mathbb{N}^n$. The \textbf{multigraded Betti number}
    \begin{equation*}
        \beta_{i,\mathbf{a}}(I) = \# \left\{ \text{Elements in a homogeneous free basis of $\mathbf{F}_i$ with multidegree $\mathbf{a}$} \right\}.
    \end{equation*}
    Let $j\in \mathbb{N}$. The \textbf{graded (or ordinary) Betti number}
    \begin{equation*}
        \beta_{i,j}(I) = \# \left\{ \text{Elements in a homogeneous free basis of $\mathbf{F}_i$ with total degree $j$} \right\} = \sum_{|\mathbf{a}|=j} \beta_{i,\mathbf{a}}(I).
    \end{equation*}
\end{defi}

\begin{example}\thlabel{resexample}
    Let $I=(x_1x_2,x_2x_3)\subset S=k[x_1,x_2,x_3]$ with grading $(1,1,1)$. Its minimal free resolution is
    \begin{equation*}
        \begin{tikzcd}[ampersand replacement=\&]
            0 \arrow[r] \& S \arrow[r, "{
                \begin{pmatrix}
                    -x_3 \\
                    x_1
                \end{pmatrix}
            }"] \& S^2 \arrow[rr, "{
                (x_1x_2 \; x_2x_3)
            }"] \&\& S \arrow[r] \& S/I \arrow[r] \& 0\\
            \substack{\text{homological} \\ \text{degree}}\& 2 \& 1 \&\& 0 \& \&
        \end{tikzcd}
    \end{equation*}
    where
    \[
        \mathbf{F}_0 = S,\,\, \mathbf{F}_1 = S[-(1,1,0)]\oplus S[-(0,1,1)],\,\, \mathbf{F}_2 = S[-(1,1,1)].
    \]
    Therefore, the nonzero Betti numbers are
    \[
        \beta_{0,(0,0,0)}(I) = 1,\,\, \beta_{1,(1,1,0)}(I) = 1,\,\,\beta_{1,(0,1,1)}(I) = 1,\,\,\beta_{2,(1,1,1)}(I) = 1,
    \]
    and
    \[
        \beta_{0,0}(I) = 1,\,\, \beta_{1,2}(I) = 2,\,\, \beta_{2,3}(I) = 1.
    \]
    If we instead use an arbitrary grading $(a,b,c)$, then the multigraded Betti numbers become
    \[
        \beta_{0,(0,0,0)}(I) = 1,\,\, \beta_{1,(a,b,0)}(I) = 1,\,\, \beta_{1,(0,b,c)}(I) = 1,\,\, \beta_{2,(a,b,c)}(I) = 1.
    \]
\end{example}
\subsection{Graphs and clique complexes}
For readers familiar with edge ideals and regularity, all the definitions in this section are standard until \thref{contribute}.

A simple graph $G$ consists of a vertex set $V(G)$ and edge set $E(G)$. We use $G^c$ to denote the complement graph of $G$. If $W$ is a subset of $V(G)$, we use $G[W]$ to denote the induced subgraph on the vertices of $W$.

\begin{defi}
    Let $G$ be a simple graph with $V(G)=\{x_1,\cdots,x_n\}$. Its \textbf{edge ideal} is
    \begin{equation*}
        I_G=\left(x_ix_j\mid \{x_i,x_j\}\in E(G)\right) \subset S=k[x_1,\cdots,x_n].
    \end{equation*}
\end{defi}

The edge ideal is a square-free quadratic monomial ideal. All square-free quadratic monomial ideals can be obtained in this way. We allow the existence of isolated vertices.

\begin{defi}
    Let $G$ be a graph with $V(G)=\{x_1,x_2,\cdots,x_n\}$. The clique complex of $G$, denoted by $\Delta G$, is a simplicial complex that consists of $(t-1)$-simplices $(x_{i_1},x_{i_2},\cdots ,x_{i_t})$ whenever the induced subgraph $G[\{x_{i_1},x_{i_2},\cdots ,x_{i_t}\}]$ is a complete graph of $t$ vertices.
\end{defi}

One can observe that the edge ideal $I_G$ is the Stanley-Reisner ideal of $\Delta G^c$. Hence, we have a one-to-one correspondence between the family of simple graphs and the family of clique complexes. Moreover, the homology of simplicial complexes $\Delta G^c$ gives Betti numbers of $I_G$ via Hochster's formula.

\begin{theorem}[Hochster's formula]\thlabel{hochster}
    Let $I_G$ be the edge ideal of a graph $G$, under standard grading. Then for $i\geq 0$ and square-free degree $\mathbf{b}$,
    \begin{equation*}
        \beta_{i,\mathbf{b}}(I_G)=\dim_k(\widetilde{H}_{|\mathbf{b}|-i-1}\left(\Delta G^c[\mathrm{supp}(\mathbf{b})], k\right)),
    \end{equation*}
    and $\beta_{i,\mathbf{b}}(I_G)=0$ if $\bf{b}$ is not square-free. If we forget the multigrading, then for $j\geq0$,
    \begin{equation*}
        \beta_{i,i+j}(I_G)=\sum_{|W|=i+j}\dim_k(\widetilde{H}_{j-1}\left(\Delta G^c[W], k\right)),
    \end{equation*}
    where $W$ runs over all the subsets of $V(G)$ of size $i+j$.
\end{theorem}

\begin{remark}
    The largest possible multidegree to have a nonzero Betti number is $(1,1,\cdots,1)$. We call this the \emph{top-level (multigraded) Betti numbers}.
\end{remark}
More detailed information on Hochster's formula can be found in \cite{francisco2014survey} and \cite{stanley1996combinatorics}. By Hochster's formula, the key part of finding the regularity of an edge ideal quotient is to find the largest number $j$, such that the $(j-1)$-st homology group of an induced subcomplex of $\Delta G^c$ is non-trivial.

\begin{defi}[Castelnuovo-Mumford regularity]\thlabel{regdef}
    For a graph $G$, the Castelnuovo-Mumford regularity of the standard graded quotient $k[G]/I_G$ is
    \begin{align*}
    \mathrm{reg}(G)&=\max\{j\mid\beta_{i,i+j}(I_G)\neq 0\}\\
    &= \max\{|\mathbf{a}|-i\mid\beta_{i,\mathbf{a}}(I_G)\neq0\}.
    \end{align*}
\end{defi}

An easy but important consequence of Hochster's formula is that multigraded Betti numbers of an induced subgraph do not depend on the ambient graph, and consequently, regularity is monotone with respect to (induced) inclusion.

\begin{proposition}\thlabel{regmono}
    Let $G$ be a graph and $H\subset V(G)$, then $\beta_{t,\mathbf{1}_H}(I_G) = \beta_{t,\mathbf{1}_H}(I_{G[H]})$. As a result, we always have $\reg(G)\geqslant \reg(G[H])$.
\end{proposition}

\begin{proof}
    By Hochster's formula we have
    \begin{equation*}
        \beta_{t,\mathbf{1}_H}(I_G) = \dim_k \tilde{H}_{|H|-t-1}(\Delta G^c[H];k),
    \end{equation*}
    and the complex $\Delta G^c[H]$ only depends on the induced subgraph $G[H]$, so the first statement holds. Now we have
    \begin{equation*}
        \reg(G) =\max\{|\mathbf{a}|-i \mid \beta_{i,\mathbf{a}}(I_G)\neq 0\}.
    \end{equation*}
    The previous statement implies that all nonzero Betti numbers of $G[H]$ also appear for $G$ at the same degree, so $\reg(G[H])$ is the maximum over a subset of the corresponding set for $\reg(G)$, and so $\reg(G[H])\leqslant \reg(G)$.
\end{proof}

For a square-free multi-index $\mathbf{a}=(a_1,\cdots,a_n)$, let $\mathrm{supp}(\mathbf{a})=\{x_i\mid a_i=1\}$.
\begin{defi}\thlabel{contribute}
    Let $V(G)=\{x_1,\cdots,x_n\}$. An induced subgraph $H\subset G$ \textbf{contributes to regularity} if $V(H)= \mathrm{supp}(\mathbf{a})$ for some multi-index $\mathbf{a}$ that attains the maximum in the definition of $\mathrm{reg}(G)$ (\thref{regdef}).
\end{defi}

\begin{example}\thlabel{resexample2}
    As in Section \ref{notation}, $P_n$ is the path with $n$ vertices, with $V(P_n)=\{x_1,\cdots,x_n\}$.

    The ideal $I_{P_3}$ is the ideal of \thref{resexample}, so $\mathrm{reg}(P_3)=1$. The multi-indices that attain the maximum are $(1,1,0),(0,1,1)$ and $(1,1,1)$, so there are three subgraphs that contribute to regularity: the two copies of $P_2$ and $P_3$ itself.

    The ideal $I_{P_4}=(x_1x_2,x_2x_3,x_3x_4)\subset S=k[x_1,x_2,x_3,x_4]$. Its non-zero Betti numbers are
    \begin{align*}
        \beta_{0,(0,0,0,0)}(I)& = 1,\,\,\beta_{1,(1,1,0,0)}(I) = 1,\,\,\beta_{1,(0,1,1,0)}(I) = 1,\\\beta_{1,(0,0,1,1)}(I) &= 1,\,\,\beta_{2,(1,1,1,0)}(I) = 1,\,\, \beta_{2,(0,1,1,1)}(I) = 1,
    \end{align*}
    so $\mathrm{reg}(P_4)=1$. There are five subgraphs that contribute to regularity: the three copies of $P_2$ and two copies of $P_3$.
\end{example}

\subsection{The lcm-lattice}

The lcm-lattice, introduced in \cite{welker1999lcm}, is a convenient way to compute the Betti numbers and analyze the minimal free resolution of edge ideals.

\begin{defi}[\cite{welker1999lcm}, Page 2]
    Let $S=k[x_1,\cdots,x_n]$ be the polynomial ring over a field $k$ and $I=(m_1,\cdots, m_d)$ a monomial ideal with its minimal generators. The elements of the \textbf{lcm-lattice} $L_I$ are the lcm's of all of the subsets of $\{m_1,... ,m_d\}$, ordered by divisibility. Its least element is $1$, which is regarded as the lcm of the empty set.
\end{defi}

Recall that an atom of a lattice is an element $x$ such that $0<x$ and there exists no element $y$ with $0<y<x$. It follows immediately that the atoms of $L_I$ are precisely $\{m_1,\cdots,m_d\}$.

The following relabeling construction from \cite{welker1999lcm} captures the adding variables construction that will be discussed in the next section.

\begin{construction}\thlabel{construction}

Let $I\subset S=k[x_1,\cdots,x_n]=k[\mathbf{x}]$ and $J\subset T=k[y_1,\cdots,y_m]=k[\mathbf{y}]$ be monomial ideals, and $f:L_I\to L_J$ a map which is a bijection on atoms and preserves joins. Let $\mathbf{P}$ be a bounded chain complex of free $S$-modules with homogeneous differentials $\partial^\mathbf{P}$, and the multidegree of every generator of every $\mathbf{P}_n$ lies in $L_I$. We construct a chain complex $f(\mathbf{P})$ of free $T$-modules as follows:

Let $\mathbf{P}_n=\bigoplus_{i=1}^{d} S[-\mathbf{a}_i]^{k_i}$, then $f(\mathbf{P})_n=\bigoplus_{i=1}^{d} T[-f(\mathbf{a}_i)]^{k_i}$

The differential $\partial^\mathbf{P}_n$ is multihomogeneous, and thus is a matrix of monomials in $S$. Suppose that its $(i,j)$-th entry $c \mathbf{x}^\mathbf{t}$ is nonzero, with $c\in k$. Suppose the $i$-th summand of $\mathbf{P}_{n-1}$ has generator in degree $\mathbf{a}$, and the $j$-th summand of $\mathbf{P}_n$ has generator in degree $\mathbf{b}$, then homogeneity implies $\mathbf{t} = \mathbf{b} - \mathbf{a}$. The differential $\partial^{f(\mathbf{P})}_n$ is the matrix having $(i,j)$-th entry $c \mathbf{y}^{f(\mathbf{b}) - f(\mathbf{a})}$. This makes it the unique multihomogeneous differential of $f(\mathbf{P})$ having the same $k$-coefficients as $\partial^\mathbf{P}_n$ in every entry.

\end{construction}

The following theorem from \cite{welker1999lcm} states that ideals with isomorphic lcm-lattice have similar minimal free resolutions.

\begin{theorem}[\cite{welker1999lcm}, Theorem 3.3]\thlabel{relabel}
    Let $I\subset S$ and $J\subset T$ be monomial ideals. Suppose that $f:L_I\to L_J$ is a bijection of the underlying sets and preserves joins. Let $\mathbf{F}$ be the minimal free resolution of $S/I$, then the multidegree of every generator of every $\mathbf{F}_n$ lies in $L_I$. Therefore, the construction $f(\mathbf{F})$ is possible. Moreover, $f(\mathbf{F})$ is precisely the minimal free resolution of $T/J$.
\end{theorem}

\begin{example}
    (Continuing \thref{resexample}). The minimal free resolution of $I=(x_1x_2,x_2x_3)\subset S=k[x_1,x_2,x_3]$ is
    \begin{equation*}
        \begin{tikzcd}[ampersand replacement=\&]
            \mathbf{F}= \& 0 \arrow[r] \& S \arrow[r, "{
                \begin{pmatrix}
                    -x_3 \\
                    x_1
                \end{pmatrix}
            }"] \& S^2 \arrow[rr, "{
                (x_1x_2 \; x_2x_3)
            }"] \&\& S \arrow[r] \& S/I \arrow[r] \& 0.\\
            \& \substack{\text{homological} \\ \text{degree}}\& 2 \& 1 \&\& 0 \& \&
        \end{tikzcd}
    \end{equation*}
    The lcm-lattice $L_I$ is (with $1$ omitted)
    \begin{center}

\tikzset{every picture/.style={line width=0.75pt}} 

\begin{tikzpicture}[x=0.75pt,y=0.75pt,yscale=-1,xscale=1]

\draw    (150,79.48) -- (179.92,109.52) ;
\draw    (150,79.48) -- (120.32,109.92) ;
\draw  [fill={rgb, 255:red, 0; green, 0; blue, 0 }  ,fill opacity=1 ] (148.32,79.48) .. controls (148.32,78.55) and (149.07,77.8) .. (150,77.8) .. controls (150.93,77.8) and (151.68,78.55) .. (151.68,79.48) .. controls (151.68,80.41) and (150.93,81.16) .. (150,81.16) .. controls (149.07,81.16) and (148.32,80.41) .. (148.32,79.48) -- cycle ;
\draw  [fill={rgb, 255:red, 0; green, 0; blue, 0 }  ,fill opacity=1 ] (118.64,109.92) .. controls (118.64,108.99) and (119.39,108.24) .. (120.32,108.24) .. controls (121.25,108.24) and (122,108.99) .. (122,109.92) .. controls (122,110.85) and (121.25,111.6) .. (120.32,111.6) .. controls (119.39,111.6) and (118.64,110.85) .. (118.64,109.92) -- cycle ;
\draw  [fill={rgb, 255:red, 0; green, 0; blue, 0 }  ,fill opacity=1 ] (178.24,109.52) .. controls (178.24,108.59) and (178.99,107.84) .. (179.92,107.84) .. controls (180.85,107.84) and (181.6,108.59) .. (181.6,109.52) .. controls (181.6,110.45) and (180.85,111.2) .. (179.92,111.2) .. controls (178.99,111.2) and (178.24,110.45) .. (178.24,109.52) -- cycle ;

\draw (126.53,63.09) node [anchor=north west][inner sep=0.75pt]    {$x_{1} x_{2} x_{3}$};
\draw (102.53,116.43) node [anchor=north west][inner sep=0.75pt]    {$x_{1} x_{2}$};
\draw (165.2,116.43) node [anchor=north west][inner sep=0.75pt]    {$x_{2} x_{3}$};

\end{tikzpicture}

    \end{center}
    Take any monomial ideal $J=(m_1,m_2)\subset T=k[y_1,\cdots,y_t]$ with $m_1\nmid m_2$ and $m_2\nmid m_1$, then $L_I$ is isomorphic to $L_J$ via the map $f:x_1x_2\mapsto m_1,\,x_2x_3\mapsto m_2$. The chain complex $f(\mathbf{F})$ below is the minimal free resolution of $T/J$.
    \begin{equation*}
        \begin{tikzcd}[ampersand replacement=\&,column sep=small]
            f(\mathbf{F})= \& 0 \arrow[r] \& T \arrow[rrr, "{
                \begin{pmatrix}
                    -\mathrm{lcm}(m_1,m_2)/m_1 \\
                    \mathrm{lcm}(m_1,m_2)/m_2
                \end{pmatrix}
            }"] \&\&\& T^2 \arrow[rr, "{
                (m_1 \; m_2)
            }"] \&\& T \arrow[r] \& T/J \arrow[r] \& 0.\\
            \& \substack{\text{homological} \\ \text{degree}}\& 2 \&\&\& 1 \&\& 0 \& \&
        \end{tikzcd}
    \end{equation*}
\end{example}

\begin{example}
Let $I=(m_1,\cdots,m_t)\subset R=k[x_1,\cdots,x_n]$. Suppose that $m_1,\cdots,m_t$ are pairwise coprime monomials in $R$, so the lcm-lattice $L_I$ is the Boolean lattice on $t$ atoms. Let $J=(y_1,\cdots,y_t)\subset S=k[y_1,\cdots,y_t]$, then $L_J$ is also the Boolean lattice on $t$ atoms, isomorphic to $L_I$ via $f:y_i\mapsto m_i$. The minimal resolution of $J$ is the Koszul complex $\mathbf{K}(y_1,\cdots,y_t)$. The construction $f(\mathbf{K}(y_1,\cdots,y_t))$ is precisely the Koszul complex $\mathbf{K}(m_1,\cdots,m_t)$ on $m_1,\cdots,m_t$, which is indeed the minimal free resolution of $I$.
\end{example}

\subsection{Splittable ideals} \label{splittable}

The material in this section will only be used in Section \ref{forests}. Splittable ideals were first introduced in \cite{ELIAHOU19901}. When a splitting $I=J+K$ exists, the Betti numbers of $I$ can be obtained from the Betti numbers of $J$, $K$ and $J\cap K$, thus allowing recursive computations of Betti numbers.

\begin{defi}
    Let $\mathcal{G}(I)$ be the set of minimal generators of a monomial ideal $I$. A monomial ideal $I$ is \emph{splittable} if there exist monomial ideals $J$ and $K$ such that
    \begin{enumerate}
        \item $I=J+K$;
        \item $\mathcal{G}(I) = \mathcal{G}(J) \sqcup \mathcal{G}(K)$, where $\sqcup$ denotes disjoint union;
        \item there exists a splitting function \begin{align*}
            \mathcal{G}(J\cap K) &\to \mathcal{G}(J) \times \mathcal{G}(K)\\
            w &\mapsto (\phi(w),\psi(w))
        \end{align*}
        such that
        \begin{enumerate}
            \item for all $w\in \mathcal{G}(J\cap K)$, $w=\lcm (\phi(w),\psi(w))$;
            \item for every nonempty subset $S\subseteq \mathcal{G}(J\cap K)$, both $\lcm (\phi(S))$ and $\lcm (\psi(S))$ strictly divide $\lcm(S)$.
        \end{enumerate}
    \end{enumerate}
    We say that $J+K$ is a \emph{splitting} of $I$ in the above case.
\end{defi}

The following theorem was proven in \cite{ELIAHOU19901} and \cite{FATABBI200192} for ungraded and $\mathbb{Z}$-graded Betti numbers. The proof in \cite{FATABBI200192} also works for multigraded Betti numbers, so we present the following theorem from \cite{FATABBI200192} in multigraded format.

\begin{theorem}[\cite{FATABBI200192}, Proposition 3.2] \thlabel{splitpre}
    Let $I=J+K$ be a splitting, then for all $i\geqslant 2$ and multidegree $\mathbf{a}$,
    \begin{align*}
        \beta_{i,\mathbf{a}}(I) &= \beta_{i,\mathbf{a}}(J) + \beta_{i,\mathbf{a}}(K) + \beta_{i-1,\mathbf{a}}(J\cap K),\\
        \beta_{1,\mathbf{a}}(I) &= \beta_{1,\mathbf{a}}(J) + \beta_{1,\mathbf{a}}(K).
    \end{align*}
\end{theorem}

Applying the above results to graphs yields the following.

\begin{theorem}[\cite{split}, Theorem 3.4 and Corollary 3.3] \thlabel{splitthm}
    Let $G$ be a graph and $e=x_ix_j$ an edge, then the following are equivalent.
    \begin{itemize}
        \item $I_G = (x_ix_j) + I_{G\setminus e}$ is a splitting.
        \item $N(x_i)\subseteq N(x_j)\cup \{x_j\}$ or $N(x_j)\subseteq N(x_i)\cup \{x_i\}$, where $N(\cdot)$ is the open neighborhood.
    \end{itemize}
    Moreover, let $H=G\setminus(N(x_i)\cup N(x_j))$, then
    \begin{equation*}
        (x_ix_j)\cap I_{G\setminus e} = x_ix_j(\, x_t\mid t\in N(x_i)\cup N(x_j)\setminus \{x_i,x_j\}\,) + x_ix_jI_H.
    \end{equation*}
\end{theorem}

\section{Regularity Offset of a Graph}
\subsection{Adding variables and adding degrees} \label{avad}


Let $G$ be a graph with $n$ vertices and $S=k[x_1,\cdots,x_n]$. Let $I_G$ be its ordinary edge ideal and $K_G$ be its generalized edge ideal with nonconstant monomials $m_1,\cdots,m_n\in T=k[z_1,\cdots,z_t]$ assigned to each vertex of $G$. Let $k_i=\deg m_i$ and let $\mathbf{g}=(k_1,\cdots,k_n)$. As the monomials $m_i$ consist of disjoint sets of variables, the lcm-lattices $L_{I_G}$ and $L_{K_G}$ are isomorphic via the bijection $f(x_ix_j)=m_im_j$.  Observe that $|f(\mathbf{a})|=\langle\mathbf{a}, \mathbf{g}\rangle$. Therefore, by \thref{relabel}, the corresponding free modules in the free resolution of $S/I_G$ and the free resolution of $T/K_G$ have the same number of generators, and a generator of degree $\mathbf{a}$ in the former corresponds to a generator of degree $f(\mathbf{a})$ in the latter. In other words, we have

\begin{proposition}\thlabel{fres}
    With the notation in the preceding paragraph, the multigraded Betti numbers satisfy
    \begin{equation*}
        \beta_{i,\mathbf{a}}(S/I_G) = \beta_{i,f(\mathbf{a})}(T/K_G).
    \end{equation*}
\end{proposition}

Now consider a new grading on $k[x_1,\dots,x_n]$ with $\deg x_i=k_i$. We have that $\deg \mathbf{x}^\mathbf{a}=\mathbf{a}\cdot \mathbf{g}=|f(\mathbf{a})|.$ Therefore, the degree of every generator of every free module in the resolution of $S/I_G$, under the grading $\deg x_i=k_i$, is the same as the degree of the corresponding generator in the resolution of $T/K_G$, under the grading $\deg z_i=1$. Thus the $\mathbb{Z}$-graded Betti numbers of the generalized edge ideal are the same as $\mathbb{Z}$-graded Betti numbers of the ordinary edge ideal, under the grading $\deg x_i=k_i$. In particular, $\mathrm{reg}(T/K_G)=\mathrm{reg}(G,\mathbf{g})$. From now on, we will only consider ordinary edge ideals under arbitrary $\mathbb{Z}$-gradings, denoting the weighted ideal by $J_G$ as in Section \ref{notation}.

Let $\mathbf{g}=(k_1,\cdots,k_n)$ be a grading. Then the degree of $x_i$ increases by $k_i-1$ when switching from the standard grading to $\mathbf{g}$. Let
\begin{equation*}
    N=\sum_{i=1}^n(k_i-1)
\end{equation*}
be the total increment of the degrees of all variables.

By Hochster's formula (\thref{hochster}), all non-zero Betti numbers of $S/I_G$ lie in square-free multidegrees. For a square-free multidegree $\mathbf{a}$, we obviously have $|\mathbf{a}|\leqslant|f(\mathbf{a})| \leqslant |\mathbf{a}| + N$. That is, after increasing the total degree by $N$, the $\mathbb{Z}$-degree of every generator in the minimal free resolution of $S/I_G$ increases by at most $N$. Thus we have shown:

\begin{theorem}\thlabel{plusn}
    Let $I_G$ be the ordinary edge ideal of a graph $G$ and $J_G$ be its edge ideal under the grading $\mathbf{g}$, whose total variable degree is $N$ more than the standard grading. Then
    \begin{equation*}
        \mathrm{reg}(G) \leqslant \mathrm{reg}(G,\mathbf{g}) \leqslant \mathrm{reg}(G)+N.
    \end{equation*}
\end{theorem}

\begin{example}\thlabel{resexample3}
    [Continuing with \thref{resexample2}.]

    For $P_3$ with grading $\mathbf{g}=(k_1,k_2,k_3)$, the nonzero Betti numbers are
    \begin{equation*}
        \beta_{0,(0,0,0)},\;\beta_{1,(k_1,k_2,0)},\;\beta_{1,(0,k_2,k_3)},\;\beta_{2,(k_1,k_2,k_3)}.\;
    \end{equation*}
    Therefore
    \begin{align*}
        \mathrm{reg}(P_3,\mathbf{g})&=\max\{k_1+k_2-1,k_2+k_3-1,k_1+k_2+k_3-2\}\\
        &= k_1+k_2+k_3-2\\
        &= \mathrm{reg}(P_3) + N.
    \end{align*}

    For $P_4$ with grading $\mathbf{g}=(k_1,k_2,k_3,k_4)$, the nonzero Betti numbers are
    \begin{align*}
        \beta_{0,(0,0,0,0)},\;\beta_{1,(k_1,k_2,0,0)},\;\beta_{1,(0,k_2,k_3,0)},\;\beta_{1,(0,0,k_3,k_4)},\;\beta_{2,(k_1,k_2,k_3,0)},\;\beta_{2,(0,k_2,k_3,k_4)}.
    \end{align*}
    Therefore
    \begin{align*}
        \mathrm{reg}(P_4,\mathbf{g})&=\max\{k_1+k_2+k_3-2,k_2+k_3+k_4-2\}\\
        &=\mathrm{reg}(P_4)+\max\{k_1+k_2+k_3-3,k_2+k_3+k_4-3\}\\
        &\leqslant\mathrm{reg}(P_4)+N,
    \end{align*}
    as $N=k_1+k_2+k_3+k_4-4$.
\end{example}

\subsection{Regularity offset of a graph} \label{tro}
The regularity grows when coordinatewise degrees increase, and \thref{plusn} provides an upper bound. We introduce the regularity offset to evaluate how far the regularity is from the upper bound.

\begin{defi}[Regularity offset]\thlabel{ro}
    Let $G$ be a graph with $n$ vertices, and $\mathbf{g}=(k_1,\cdots,k_n)\in(\mathbb{N}^+)^n$ a grading, and $N=\sum_{i=1}^{n} (k_i-1)$. Let $I_G$ be the ordinary edge ideal and $J_G$ be the edge ideal with grading $\mathbf{g}$. We define the regularity offset of $G$ under grading $\mathbf{g}$ as
    \begin{equation*}
        \mathrm{ro}(G,\mathbf{g})=\mathrm{reg}(G)+N-\mathrm{reg}(G,\mathbf{g})
    \end{equation*}
    We say $G$ does not have regularity offset if $\mathrm{ro}(G,\mathbf{g})= 0$ for all gradings $\mathbf{g}\in(\mathbb{N}^+)^n$. In other words, $\mathrm{reg}(G,\mathbf{g})$ always attains maximum growth.
\end{defi}

\begin{example}
    [Continuing with \thref{resexample3}.] Let $t_i=k_i-1$ be the increment of the degree of $x_i$. Then we have
    \begin{align*}
        \mathrm{ro}(P_3,\mathbf{g}) = 0,\;\mathrm{ro}(P_4,\mathbf{g}) = \min\{t_1,t_4\}.
    \end{align*}
\end{example}

To keep track of all Betti supports, write $r=\mathrm{reg}(G)$ and let
\[
    \mathcal{B}(G)=\{(i,W):i\geq0,\ W\subseteq V(G),\ \beta_{i,\mathbf{1}_W}(I_G)\neq0\}.
\]
This set includes $(0,\varnothing)$. For a grading $\mathbf{g}$, put $t_v=\deg_{\mathbf{g}}v-1$ and $N=\sum_{v\in V(G)}t_v$. Regrading changes the total degree of a square-free Betti support $W$ from $|W|$ to $|W|+\sum_{v\in W}t_v$, without changing its homological degree. Hence
\begin{equation}\label{eq:weighted-supports}
    \mathrm{reg}(G,\mathbf{g})
    =\max_{(i,W)\in\mathcal{B}(G)}
      \left(|W|-i+\sum_{v\in W}t_v\right).
\end{equation}
Subtracting this maximum from $r+N$ gives
\begin{equation}\label{eq:offset-supports}
    \mathrm{ro}(G,\mathbf{g})
    =\min_{(i,W)\in\mathcal{B}(G)}
      \left(r-(|W|-i)+\sum_{v\notin W}t_v\right).
\end{equation}
Both $r-(|W|-i)$ and the sum in each term are nonnegative. This formula shows that maximal growth for the single grading $(2,\ldots,2)$ already implies maximal growth for every positive grading.

\begin{lemma}\thlabel{rogt0}
    A graph $G$ with $n$ vertices does not have regularity offset if and only if
    \begin{equation*}
        \mathrm{ro}(G,\mathbf{t})=0,
    \end{equation*}
    where $\mathbf{t}\in\mathbb{Z}^n$, $\mathbf{t}=(2,\cdots,2)$.
\end{lemma}

\begin{proof}
    If $G$ has no regularity offset, then $\mathrm{ro}(G,\mathbf{t})=0$ by definition. Conversely, suppose that $\mathrm{ro}(G,\mathbf{t})=0$. For this grading every increment $t_v$ is $1$, so \eqref{eq:offset-supports} gives
    \[
        0=\min_{(i,W)\in\mathcal{B}(G)}
          \bigl(r-(|W|-i)+|V(G)\setminus W|\bigr).
    \]
    Since both summands are nonnegative, a pair attaining this minimum must satisfy $W=V(G)$ and $|W|-i=r$. For any positive grading $\mathbf{g}$, this same pair contributes zero to the minimum in \eqref{eq:offset-supports}. Thus $\mathrm{ro}(G,\mathbf{g})=0$ for every such grading.
\end{proof}

With this lemma, we give a topological condition for existence of the regularity offset.

\begin{theorem}\thlabel{wholegraph}
    A graph $G$ does not have regularity offset if and only if the whole graph contributes to regularity, or equivalently, $\widetilde{H}_{\mathrm{reg}(G)-1}\left(\Delta G^c, k\right)$ is non-trivial.
\end{theorem}

\begin{proof}
    The proof of \thref{rogt0} shows that $G$ has no regularity offset if and only if there exists $(i,V(G))\in\mathcal{B}(G)$ with $|V(G)|-i=\mathrm{reg}(G)$. This is exactly the assertion that the whole graph contributes to regularity. Writing $n=|V(G)|$ and $r=\mathrm{reg}(G)$, Hochster's formula gives
    \[
        \beta_{n-r,\mathbf{1}_{V(G)}}(I_G)
        =\dim_k\widetilde{H}_{r-1}(\Delta G^c;k),
    \]
    proving the homological characterization.
\end{proof}

After understanding existence of regularity offset, we are interested in where the regularity offset comes from, or in other words, which are the vertices that cause the offset?

\begin{defi}\thlabel{offsetvertices}
    Let $\mathcal{R}_{\max}(G)$ be the inclusion-maximal subsets $W\subseteq V(G)$ for which $\beta_{i,\mathbf{1}_W}(I_G)\neq0$ and $|W|-i=\mathrm{reg}(G)$ for some $i$. The set of \textbf{offset vertices} is
    \[
        V_o(G)=\bigcup_{W\in\mathcal{R}_{\max}(G)}\bigl(V(G)\setminus W\bigr).
    \]
    Thus an offset vertex is omitted from at least one maximal regularity-contributing support. This definition uses the standard grading; weighted offsets are computed from all pairs in $\mathcal{B}(G)$ using \eqref{eq:offset-supports}.
\end{defi}

\subsection{\texorpdfstring{The lower bound in \thref{plusn}}{The lower bound for regularity}}
\thref{rogt0} and \thref{wholegraph} gave conditions on when the upper bound in \thref{plusn} is attained. Now we show that the lower bound can also be attained for any $N$.
\begin{proposition}\thlabel{prlower}
    For any integer $N\geq0$, there exist a graph $G$ and a positive grading $\mathbf{g}=(k_1,\cdots,k_n)$ with $\sum_{i=1}^n(k_i-1)=N$ such that $\mathrm{reg}(G,\mathbf{g})=\mathrm{reg}(G)$.
\end{proposition}

\begin{proof}
We introduce a family of graphs such that, for any $N$, there exists a graph in this family whose regularity remains the same after increasing total variable degree by $N$.
  \begin{figure}[h]
      \centering
      \includegraphics[width=0.8\linewidth]{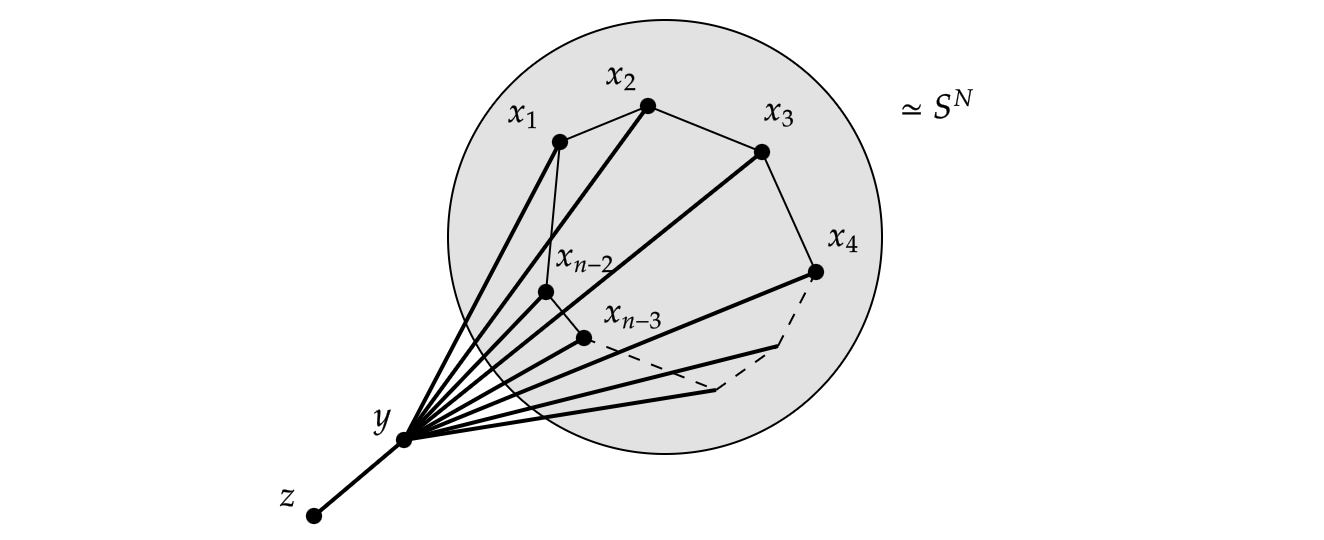}
  \end{figure}

The figure schematically illustrates the construction. Let $H$ be the disjoint union of $N+1$ edges on vertices $x_1,\ldots,x_{n-2}$, where $n=2N+4$. Its edge ideal is generated by $N+1$ quadratic monomials in disjoint variables, so its quotient has the Koszul resolution and $\mathrm{reg}(H)=N+1$. Also, $\Delta H^c$ is the join of $N+1$ copies of $\mathbb{S}^0$, hence is homeomorphic to $\mathbb{S}^N$. Form $G$ by adding vertices $y,z$, with $y$ adjacent to every other vertex and $z$ adjacent only to $y$.

We claim that $\mathrm{reg}(G)=N+1$. Since $H$ is induced in $G$, \thref{regmono} gives the lower bound. For the upper bound, consider any induced subcomplex of $\Delta G^c$. If it contains $z$, its part without $y$ is a cone with apex $z$. If it does not contain $z$, its part without $y$ is an induced subcomplex of $\Delta H^c$. Including $y$ only adds an isolated point. Thus no such subcomplex has nonzero reduced homology above dimension $N$, and Hochster's formula gives $\mathrm{reg}(G)\leq N+1$.

Now assign $z$ degree $N+1$ and every other vertex degree $1$, obtaining a grading $\mathbf{g}$ with total increment $N$. A nonzero Betti number whose support omits $z$ retains its degree and contributes at most $N+1$ to $\mathrm{reg}(G,\mathbf{g})$. If its support contains $z$, it must also contain $y$; the corresponding induced complex is a nonempty cone together with an isolated point. Its only nonzero reduced homology is in degree $0$, so its standard contribution to regularity is $1$. Regrading increases its total degree by $N$ without changing its homological degree, making its contribution $N+1$. Consequently,
\[
    \mathrm{reg}(G,\mathbf{g})=N+1=\mathrm{reg}(G).
\]
\end{proof}


\section{Regularity Offset of Some Simple Graphs} \label{rossg}

\subsection{Co-chordal graphs}
Chordal graphs are graphs with no induced cycles of length larger than $3$. We call a graph co-chordal if its complement graph is chordal. Co-chordal graphs have linear resolutions.
\begin{theorem}[\cite{RalfFröberg1990}, Fröberg’s Theorem]
    Let $G$ be a graph with at least one edge. Its edge ideal has a linear resolution if and only if $G^c$ is chordal. In particular, a co-chordal graph with at least one edge has $\mathrm{reg}(G)=1$; an edgeless graph has $\mathrm{reg}(G)=0$.
\end{theorem}

This result allows us to determine when co-chordal graphs have no regularity offset, and also compute their regularity offset.

\begin{proposition}
    If $G$ is co-chordal and has at least one vertex, then $G$ does not have regularity offset if and only if $G^c$ is not connected.
\end{proposition}

\begin{proof}
    If $G$ is edgeless and has at least one vertex, then $\mathrm{reg}(G)=\mathrm{reg}(G,\mathbf{g})=0$ for every grading. Thus $\mathrm{ro}(G,\mathbf{g})=\sum_i(k_i-1)$, so $G$ has regularity offset, and $G^c$ is connected. Otherwise, $G$ has an edge and $\mathrm{reg}(G)=1$. By \thref{wholegraph}, $G$ does not have regularity offset if and only if $\widetilde{H}_{0}\left(\Delta G^c, k\right)$ is non-trivial. We observe that $\dim_k(\widetilde{H}_{0}\left(\Delta G^c, k\right))=\#$ connected components of  $\Delta G^c-1$. Hence, $G$ does not have regularity offset if and only if $G^c$ is not connected.
\end{proof}

\begin{theorem}\thlabel{ccd}
    Let $G$ be a co-chordal graph with a positive grading $\mathbf{g}$. Put $t_v=\deg_{\mathbf{g}}v-1\geq0$ and $N=\sum_{v\in V(G)}t_v$. Define
    \[
        \mathcal{C}(G)=\{C\subseteq V(G):G^c[V(G)\setminus C]
            \text{ has at least two connected components}\},
    \]
    and let $\mathcal{C}_{\min}(G)$ be its inclusion-minimal members. If $E(G)\neq\varnothing$, then
    \[
        \mathrm{ro}(G,\mathbf{g})
          =\min_{C\in\mathcal{C}_{\min}(G)}\sum_{v\in C}t_v.
    \]
    In particular, $\mathcal{C}_{\min}(G)=\{\varnothing\}$ when $G^c$ has at least two connected components. If $E(G)=\varnothing$, then $\mathcal{C}(G)=\varnothing$ and $\mathrm{ro}(G,\mathbf{g})=N$.
\end{theorem}
\begin{proof}
    Suppose first that $G$ has an edge, so $\mathrm{reg}(G)=1$. For a nonempty support $W$, every nonzero $\beta_{i,\mathbf{1}_W}(I_G)$ satisfies $|W|-i=1$. Indeed, the linear resolution of the edge ideal gives this equality in the quotient resolution. By Hochster's formula,
    \[
        \beta_{|W|-1,\mathbf{1}_W}(I_G)
          =\dim_k\widetilde{H}_0(\Delta G^c[W];k).
    \]
    The clique complex $\Delta G^c[W]$ has the same connected components as $G^c[W]$. Thus $W$ is a nonempty Betti support precisely when $C=V(G)\setminus W$ belongs to $\mathcal{C}(G)$. This family is nonempty: the endpoints of any edge of $G$ induce two isolated vertices in $G^c$.

    The pair $(0,\varnothing)$ contributes $1+N$ to \eqref{eq:offset-supports}, whereas each nonempty Betti support contributes $\sum_{v\in C}t_v\leq N$. Consequently,
    \[
        \mathrm{ro}(G,\mathbf{g})
          =\min_{C\in\mathcal{C}(G)}\sum_{v\in C}t_v.
    \]
    Every member of the finite family $\mathcal{C}(G)$ contains an inclusion-minimal member. Since $t_v\geq0$, restricting the minimum to these members does not change it. If $G^c$ already has at least two connected components, $\varnothing\in\mathcal{C}(G)$ and is its unique inclusion-minimal member.

    Finally, if $G$ is edgeless, both quotient regularities are zero, so the offset is $N$. Every induced subgraph of $G^c$ is complete and has at most one connected component, hence $\mathcal{C}(G)=\varnothing$.
\end{proof}
\begin{cor}
    For a co-chordal graph $G$ with at least one edge,
    \[
        V_o(G)=\bigcup_{C\in\mathcal{C}_{\min}(G)}C.
    \]
    If $G$ is edgeless, then $V_o(G)=V(G)$.
\end{cor}
\begin{proof}
    In the case with an edge, the preceding proof identifies regularity-contributing supports with the complements of members of $\mathcal{C}(G)$. Maximal such supports therefore correspond to inclusion-minimal members of $\mathcal{C}(G)$, giving the formula by \thref{offsetvertices}. In the edgeless case, the only Betti support is $\varnothing$, and the same definition gives $V_o(G)=V(G)$.
\end{proof}

\subsection{Paths and cycles}
Paths and cycles are two other families of graphs that have relatively simple structure. We discuss the behavior of regularity offset of them in this section.

First of all, we need to understand the regularity of paths and cycles.

\begin{proposition}[\cite{woodroofe2014matchings}, Proposition 10]

For $n\geq3$, $\mathrm{reg}(C_n)=\mathrm{reg}(P_n)=\left\lfloor(n+1)/3\right\rfloor$. The path formula also holds for $n=1,2$.

\end{proposition}
We discuss cycles first. By \thref{wholegraph}, to understand cycles' regularity offset, we just need to figure out the homotopy type of $\Delta C_n^c$. In Proposition 5.2 of \cite{kozlov1999complexes}, Kozlov showed that:
\begin{equation*}
    \Delta C_n^c\simeq\left\{\begin{matrix}
 \mathbb{S}^{k-1}\vee \mathbb{S}^{k-1},\qquad &\text{if}\quad  n=3k;\\
 \mathbb{S}^{k-1},\qquad &\text{if}\quad  n=3k+1;\\
 \mathbb{S}^k,\qquad &\text{if}\quad  n=3k+2.
\end{matrix}\right.
\end{equation*}
We can check that these homotopy types fit in the regularity exactly. Hence, we can give the conclusion for cycles.

\begin{proposition}\thlabel{cycle}
    Cycle graphs do not have regularity offset.
\end{proposition}
Now, let's move on to the paths. Similarly, we need to introduce a topological conclusion at the beginning.

\begin{proposition}[\cite{kozlov1999complexes}, Proposition 4.6]\thlabel{P}
\begin{equation*}
    \Delta P_n^c\simeq \left\{\begin{matrix}
 \mathbb{S}^{k-1},\qquad &\text{if}\quad  n=3k;\\
 \text{a point,}\qquad &\text{if}\quad  n=3k+1;\\
 \mathbb{S}^k,\qquad &\text{if}\quad  n=3k+2.
\end{matrix}\right.
\end{equation*}
\end{proposition}
 Therefore, we have the following corollary immediately.
 \begin{cor}\thlabel{path}
     $P_n$ does not have regularity offset if and only if $n\equiv0,2$ $(\mathrm{mod}$ $3)$.
 \end{cor}
For the case that $n\equiv1$ $(\mathrm{mod}$ $3)$, we have the following theorem.
\begin{theorem}\thlabel{path2}
    Let $n=3k+1$ with $k\geq0$, and label the vertices of $P_n$ in order as $x_1,\ldots,x_n$. For a positive grading $\mathbf{g}$, put $t_i=\deg_{\mathbf{g}}x_i-1$. Then
    \[
        V_o(P_n)=\{x_{3s+1}:0\leq s\leq k\},\qquad
        \mathrm{ro}(P_n,\mathbf{g})=\min_{0\leq s\leq k}t_{3s+1}.
    \]
\end{theorem}
\begin{proof}
    If $k=0$, then $P_n=P_1$ is edgeless, its only Betti support is $\varnothing$, and the assertions follow directly. Assume $k\geq1$. Write $V=V(P_n)$, $A=\{x_{3s+1}:0\leq s\leq k\}$, and $S_j=V\setminus\{x_j\}$ for $j\equiv1\pmod3$. We have $\mathrm{reg}(P_n)=k$.

    For $j=3s+1$, the components of $P_n[S_j]$ have $3s$ and $3(k-s)$ vertices, with an empty component omitted. The independence complex of a disjoint union is the join of the independence complexes of its components (see \cite{kozlov2008combinatorial}, Corollary 20.6). By \thref{P}, this join is homotopy equivalent to $\mathbb{S}^{k-1}$; when both components are nonempty, it is
    \[
        \mathbb{S}^{s-1}*\mathbb{S}^{k-s-1}\simeq\mathbb{S}^{k-1}.
    \]
    Hochster's formula therefore gives $\beta_{2k,\mathbf{1}_{S_j}}(I_{P_n})\neq0$, and $|S_j|-2k=k$. Thus every $S_j$ contributes to regularity.

    We next show that every nonzero Betti support omits a vertex of $A$. The empty support does so. For a nonempty support $W$, the components of $P_n[W]$ are paths. If one component had a number of vertices congruent to $1$ modulo $3$, \thref{P} would make its independence complex contractible. The join for $P_n[W]$ would then also be contractible, contradicting Hochster's formula. Thus each component of $P_n[W]$ has order congruent to $0$ or $2$ modulo $3$.

    Suppose, for a contradiction, that $A\subseteq W$. Write the components from left to right as intervals $[a_1,b_1],\ldots,[a_r,b_r]$ of vertex indices. Then $a_1=1$ and $b_r=n$. We prove inductively that $a_i\equiv1\pmod3$. If this holds for $i<r$, the component-order condition gives $b_i\equiv0$ or $2\pmod3$. The first case would put $x_{b_i+1}\in A\subseteq W$, contradicting the end of the component. Hence $b_i\equiv2\pmod3$. The next vertex $x_{b_i+1}$ is absent, while $x_{b_i+2}\in A\subseteq W$, so $a_{i+1}=b_i+2\equiv1\pmod3$. The last component consequently has order $n-a_r+1\equiv1\pmod3$, a contradiction.

    Every Betti support is therefore contained in some $S_j$ with $x_j\in A$. Since each $S_j$ contributes to regularity and the full vertex set is not a Betti support, these are exactly the maximal regularity-contributing supports. The assertion about $V_o(P_n)$ follows from \thref{offsetvertices}.

    Finally, for every $(i,W)\in\mathcal{B}(P_n)$, we have $k-(|W|-i)\geq0$ and $x_j\notin W$ for some $x_j\in A$. Thus
    \[
        k-(|W|-i)+\sum_{x_l\notin W}t_l
          \geq t_j\geq\min_{0\leq s\leq k}t_{3s+1}.
    \]
    Conversely, the pair $(2k,S_j)$ contributes exactly $t_j$ to \eqref{eq:offset-supports}. Taking the minimum over all pairs proves the formula for every positive grading.
\end{proof}

\subsection{Forests} \label{forests}

We will use the theory of splittable edge ideals to analyze the regularity offset of forests. Splittable ideals were introduced in \cite{ELIAHOU19901} and \cite{FATABBI200192}, along with some relevant theorems. The application of splittable ideals to edge ideals comes from \cite{split}. We collected these results in Section \ref{splittable}.

As non-zero Betti numbers occur only in square-free multidegrees, and square-free multidegrees are in one-to-one correspondence with subsets of the vertex set, we introduce the following notation. Let $G$ be a graph and let $S$ be a subset of vertices of $G$, then
\begin{equation*}
    \beta_t(G,S) \coloneq \beta_{t,\mathbf{1}_S}(I_G).
\end{equation*}

Forests with at least one edge always have a leaf, and the edge of a leaf always satisfies the second condition in \thref{splitthm}, and therefore provides a splitting of the edge ideal. In view of \thref{regmono}, to compute the multigraded Betti numbers of a graph, we only need to compute the top-level multigraded Betti numbers $\beta_t(G,V(G))$ for every induced subgraph. For forests, we have the following result.

\begin{proposition} \thlabel{deletion}
    Let $T$ be a forest. Let $x_n\in V(T)$ be a leaf and $x_{n-1}$ its neighbor. Let $S\subset T$ be the induced subforest on the vertex set $V(T)\setminus (N(x_{n-1})\cup \{x_{n-1}\})$, then
    \begin{equation*}
        \beta_t(T,V(T)) = \beta_{t-|N(x_{n-1})|}(S,V(S)).
    \end{equation*}
\end{proposition}

\begin{proof}
    Write $R=k[V(T)]$. Let $e=x_{n-1}x_n$, then by \thref{splitthm} $I_T=(x_{n-1}x_n)+I_{T\setminus e}$ is a splitting. Write $L=(x_{n-1}x_n)\cap I_{T\setminus e}$. By \thref{splitpre} we have
    \begin{equation*}
        \beta_t(T,V(T)) = \beta_{t,\mathbf{1}_{V(T)}}((x_{n-1}x_{n})) + \beta_{t,\mathbf{1}_{V(T)}}(I_{T\setminus e}) + \beta_{t-1,\mathbf{1}_{V(T)}}(L)\cdot 1_{[t\geqslant 2]},
    \end{equation*}
    where $1_{[t\geqslant 2]}=1$ when $t\geqslant 2$ and $=0$ otherwise. The minimal free resolution of $R/(x_{n-1}x_n)$ shows that the first term is $1$ when $t=1$ and $V(T)=\{x_{n-1},x_n\}$, and $0$ otherwise. As the ideal $I_{T\setminus e}$ does not contain the variable $x_n$ in any of its generators, but the entry corresponding to $x_n$ is nonzero in the multidegree $\mathbf{1}_{V(T)}$, the second term $\beta_{t,\mathbf{1}_{V(T)}}(I_{T\setminus e})$ is $0$.

    We now focus on the third term. By \thref{splitthm} we have
    \begin{equation*}
        L=x_{n-1}x_n(\,(x_i\mid x_i \in N(x_{n-1})\setminus\{x_n\}) + I_S\,).
    \end{equation*}
    Set
    \begin{align*}
        J&=(x_i\mid x_i \in N(x_{n-1})\setminus\{x_n\})\\
        M&=J + I_S.
    \end{align*}
    Then no generator of $M$ contains $x_{n-1}$ or $x_n$. Let $\mathbf{j}=\mathbf{1}_{V(T)\setminus \{x_{n-1},x_n\}}$, then
    \begin{equation*}
        \beta_{t-1,\mathbf{1}_{V(T)}}(L) = \beta_{t-1,\mathbf{j}}(M).
    \end{equation*}
    Now $J$ and $I_S$ have disjoint and pairwise coprime sets of generators. Let $R'=k[N(x_{n-1})\setminus \{x_n\}]$ and $R''=k[V(S)]$, then $R/M\cong R'/J \otimes R''/I_S \otimes k[x_{n-1},x_n]$. Therefore
    \begin{equation*}
        \beta_{t-1,\mathbf{j}}(R/M) = \sum_{a+b=t-1}\beta_{a,\mathbf{1}_{R'}}(R'/J)\beta_{b,\mathbf{1}_{R''}}(R''/I_S).
    \end{equation*}
    The minimal free resolution of $R'/J$ is the Koszul complex, so the only nonzero $\beta_{a,\mathbf{1}_{R'}}(R'/J)$ occurs at $a=|\mathrm{supp}(\mathbf{1}_{R'})|=|N(x_{n-1})|-1$, and its value is $1$. Therefore
    \begin{equation*}
        \beta_{t-1,\mathbf{j}}(R/M) = \beta_{t-|N(x_{n-1})|,\mathbf{1}_{R''}}(R''/I_S) = \beta_{t-|N(x_{n-1})|}(S,V(S)).
    \end{equation*}

    Putting together all the cases, we have
    \begin{equation*}
        \beta_t(T,V(T)) = \beta_{t-|N(x_{n-1})|}(S,V(S))
    \end{equation*}
    when $t\geqslant 2$, and
    \begin{equation*}
        \beta_t(T,V(T)) = 1_{[V(T)=\{x_{n-1},x_n\}]}
    \end{equation*}
    when $t=1$.
    When $t=1$, $\beta_{t-|N(x_{n-1})|}(S,V(S))=\beta_{1-|N(x_{n-1})|}(S,V(S))$. As $|N(x_{n-1})|\geqslant 1$, this can only be nonzero when $|N(x_{n-1})|=1$ and $S=\varnothing$. This finishes the proof.
\end{proof}

\thref{deletion} provides a completely graph-theoretic certificate for when a forest's top-level Betti number is $1$. We call these forests \emph{full forests}.

\begin{proposition} \thlabel{deletion2}
    Let $T$ be a forest. Define a deletion operation as follows: take a leaf with its neighbor, and remove all vertices in the neighborhood of its neighbor, plus the neighbor itself. Then $T$ is full if and only if it can be reduced to the empty forest after a finite number of deletions.
    If $T$ is not full, then any order of deletion operations always ends up with a set of isolated vertices.
\end{proposition}

\begin{proof}
    According to \thref{deletion}, the top-level Betti number stays constant through a deletion, with its homological degree dropping by $|N(x_{n-1})|$. The deletion operation can be performed on any forest with edges, so any forest can be reduced to either the empty graph or an edgeless set of isolated vertices after a finite number of deletions. If it can be reduced to the empty graph, then as $\beta_0(\varnothing,\varnothing)=1$, there exists a unique $t$ such that $\beta_t(T,V(T))=1$. This $t$ is equal to the sum of the $|N(x_{n-1})|$ for each operation performed. On the other hand, if $T$ is reduced to a set of isolated vertices $S$, then $\beta_l(S,V(S))=\beta_{l,\mathbf{1}}(R/0)=0$ for all $l$, so $\beta_t(T,V(T))=0$ for all $t$.
\end{proof}

\begin{proposition} \thlabel{deletion3}
    If $T$ is a full forest, then there exists a unique $t$ such that $\beta_t(T,V(T))\neq 0$. Denote this number by $i_T$. We have $i_T=|V(T)|-\reg(T)$. Let $\sigma_T$ be the number of deletions used to reduce it to the empty graph; we have $\sigma_T = \mathrm{reg}(T)$
\end{proposition}

\begin{proof}
    Almost all of the claims were verified during the proof of \thref{deletion2}. The remaining assertions to prove are $\sigma_T = \mathrm{reg}(T)$ and $i_T=|V(T)|-\reg(T)$.

    Each deletion operation reduces the homological degree by $|N(x_{n-1})|$, which is $1$ fewer than the number of vertices deleted. Therefore the whole chain of deletions reduces the homological degree by $|V(T)|-\sigma_T$, and ends with a $0$. Thus $i_T=|V(T)|-\sigma_T$.

    In each step of our deletion process, choose the edge between the selected leaf $x_i$ and its neighbor $x_j$. No endpoint of a later chosen edge can be adjacent to $x_i$ or $x_j$. This is because $x_i$ is only adjacent to $x_j$ and (possibly) earlier deleted vertices, and the entire neighborhood of $x_j$ is deleted in this step. Therefore, the set of all chosen edges forms an induced matching of cardinality $\sigma_T$.

    As $T$ is full, every vertex is deleted at some point, and the deleted set of vertices form a star or isolated edge. Therefore every vertex is either an endpoint of a chosen edge or adjacent to an endpoint of a chosen edge. This implies that this matching is maximal.

    We now prove that this induced matching is maximum. Take any other induced matching $M$. If an edge $e\in M$ is in one of the stars, then it is the only edge of $M$ in that star. If $e\in M$ is not in any star, then it is the connecting edge between two stars, as all vertices are in some star. In this case, no edge in these two stars can be in $M$ because $M$ is induced. Contract all of these stars into single vertices, and we get a forest $T^-$ with $\sigma_T$ vertices. The induces matching $M$ becomes a collection of vertices and edges of $T^-$ such that neither endpoint of a selected edge is a selected vertex. Therefore, the cardinality of $M$ is no greater than the number of vertices of $T^-$, which is $\sigma_T$.
    
    Theorem 2.18 of \cite{Zheng31122004} states that the maximum induced matching number of a forest is equal to its regularity. Therefore $\reg(T)=\sigma_T$.
\end{proof}

Now we turn to the regularity offset of forests. Let $T$ be a forest with $n$ vertices and $J_T$ be its edge ideal under the grading $\mathbf{g}=(k_1,\cdots,k_n)$.
\begin{align*}
    \mathrm{reg}(T) &= \max_{S\subset T\;\text{full}}\left( |V(S)| - i_S \right),\\
    \mathrm{reg}(T,\mathbf{g}) &= \max_{S\subset T\;\text{full}}\left( \sum_{x_l\in V(S)}k_l - i_S \right),\\
    \mathrm{ro}(T,\mathbf{g}) &= \mathrm{reg}(T) -n + \min_{S\subset T\;\text{full}} \left( \sum_{x_l\notin V(S)}k_l + i_S \right).
\end{align*}

With \thref{deletion3}, we get a better understanding of the behavior of regularity offset.

\begin{theorem} \thlabel{final}
    Let $T$ be a forest on $\{x_1,\cdots,x_n\}$, $Z$ the set of full induced subforests of $T$, $\mathbf{g}=(k_1,\cdots,k_n)$, and $t_i=k_i-1$ the increment of degree of $x_i$, then
    \begin{equation*}
        \mathrm{ro}(T,\mathbf{g}) = \mathrm{reg}(T) + \min_{S\in Z\;\text{maximal}}\left( \sum_{i:\,x_i\notin V(S)} t_i - \mathrm{reg}(S) \right).
    \end{equation*}
\end{theorem}

\begin{proof}
    By \thref{deletion3}, we have
    \begin{align*}
        \mathrm{ro}(T,\mathbf{g}) &= \mathrm{reg}(T) + \min_{S\subset T\;\text{full}} \left( \sum_{x_l\notin V(S)}k_l -n+ i_S \right)\\
        &= \mathrm{reg}(T) + \min_{S\subset T\;\text{full}} \left( \sum_{x_l\notin V(S)}k_l -n+ |V(S)|-\mathrm{reg}(S) \right)\\
        &= \mathrm{reg}(T) + \min_{S\subset T\;\text{full}} \left( \sum_{x_l\notin V(S)}t_l-\mathrm{reg}(S) \right).
    \end{align*}
    When an inclusion $S_1\subsetneq S_2$ exists, we have $\sum_{x_l\notin V(S_1)}t_l\geqslant \sum_{x_l\notin V(S_2)}t_l$ and $\mathrm{reg}(S_1)\leqslant \mathrm{reg}(S_2)$. Therefore the minimum in the formula for $\mathrm{ro}(T,\mathbf{g})$ can be taken over all maximal $S\in Z$.
\end{proof}

To illustrate \thref{final}, we re-derive \thref{path2} about the regularity offset of paths using \thref{final}.

\begin{proposition}\thlabel{pathfinal}
    Let $n=3k+1$ with $k\geq0$, so the path $P_n$ is not full (see \thref{path}). Write $V(P_n)=\{x_1,\cdots,x_n\}$ so that the $x_i$ are labeled in order. Then the maximal full subforests of $P_n$ are exactly the induced subgraphs on the sets $S_{3s+1}=V(P_n)\setminus\{x_{3s+1}\}$ for $0\leq s\leq k$.
\end{proposition}

\begin{proof}
    The full induced subforests of $P_n$ correspond exactly to its nonzero Betti supports, by \thref{regmono} and \thref{deletion2}. For $k\geq1$, the proof of \thref{path2} identifies the maximal such supports as the sets $S_{3s+1}$ with $0\leq s\leq k$. For $k=0$, the only full induced subforest is the empty graph, corresponding to $S_1=\varnothing$.
\end{proof}

\begin{cor}
    Continue with the assumptions of \thref{pathfinal}. Let $\mathbf{g}=(k_1,\cdots,k_n)$ and $t_i=k_i-1$, then
    \begin{equation*}
        \mathrm{ro}(P_n,\mathbf{g}) = \min_{0\leq s\leq k}t_{3s+1}.
    \end{equation*}
\end{cor}

\begin{proof}
    We have $\mathrm{reg}(P_{3k+1})=k$, and each maximal full subforest also has regularity $k$. The formula now follows from \thref{final}.
\end{proof}

\bibliographystyle{plain}
\bibliography{reference}
\end{document}